\documentclass[10pt,twoside]{article}
\usepackage{graphicx}
\usepackage{amsmath}
\usepackage{Latex-document}
\newtheorem{theorem}{Theorem}
\newtheorem{problem}{Problem}

\def\R{{\mathbf{R}}}
\def\rn{{\mathbf{R}^n}}

\def\cS{{\mathbf{S}}}

\title{\bf The Branch Set of a Quasiregular\vskip -2mm Mapping\thanks{Supported by NSF grant DMS
9970427. I thank Mario Bonk and Alex Eremenko for their criticism
on earlier versions of this article. My warmest thanks go to
Mario Bonk, Seppo Rickman, and Dennis Sullivan for collaboration,
mentoring, and friendship.}\vskip 6mm}
\author{Juha Heinonen\vspace*{-0.5cm}\thanks{Department of Mathematics,
University of Michigan, MI 48109, USA.  E-mail:
juha@math.lsa.umich.edu }}
\date{\vspace{-8mm}}

\begin{document}

\maketitle

\thispagestyle{first}\setcounter{page}{691}

\begin{abstract}

\vskip 3mm

We discuss the issue of branching in quasiregular mapping, and in
particular the relation between branching and the problem of
finding geometric parametrizations for topological manifolds.
Other recent progress and open problems of a more function
theoretic nature are also presented.

\vskip 4.5mm

\noindent{\bf 2000 Mathematics Subject Classification:}  30C65,
57M12.

\noindent{\bf Keywords and Phrases:}  Quasiregular map,
Bi-Lipschitz map, Branch set.
\end{abstract}

\vskip 12mm

\section{Branched coverings}

\vskip -5mm \hspace{5mm}

A continuous mapping $f:X\to Y$ between topological spaces is said
to be a {\textit{branched covering}} if $f$ is an open mapping and
if for each $y\in Y$ the preimage $f^{-1}(y)$ is a discrete subset
of $X$. The {\textit{branch set}} $B_f$ of $f$ is the closed set
of points in $X$ where $f$ does not define a local homeomorphism.

Nonconstant holomorphic functions between connected Riemann surfaces are
examples of branched coverings.
From the Weierstrassian (power series) point of view this property
of holomorphic functions is
almost immediate.
It is a deeper fact, due to Riemann,
 that the same conclusion can be drawn from the mere
definition of complex differentiability, or, equivalently, from the
Cauchy-Riemann equations.
Most of this article discusses
the repercussions of this fact.

\section{Quasiregular mappings}

\vskip -5mm \hspace{5mm}

In a  1966 paper [27], Reshetnyak penned a definition for {\it
mappings of bounded distortion} or, as they are more commonly
called today, {\it quasiregular mappings}. These are nonconstant
 mappings $f:\Omega\to\rn$
 in the Sobolev space
$W^{1,n}_{loc}(\Omega;\rn)$, where
$\Omega\subset
\rn$ is a domain and $n\ge2$,
 satisfying the following requirement: there exists a constant
$K\ge1$ such that
\begin{equation}\label{qrdef1}
|f'(x)|^n\le KJ_f(x)
\end{equation}
for almost every $x\in \Omega$, where $|f'(x)|$ denotes the operator norm of
the (formal) differential matrix
$f'(x)$ with $J_f(x)$=det$f'(x)$ its Jacobian determinant.
One also speaks about
$K$-{\textit {quasiregular mappings}} if the constant in \eqref{qrdef1}
is to be emphasized.\footnote{The definition readily
extends for mappings between connected oriented
Riemannian $n$-manifolds.}

Requirement \eqref{qrdef1}  had been
used as the analytic definition for quasiconformal mappings since
the 1930s, with varying degrees of smoothness conditions
on $f$. Quasiconformal mappings are by definition quasiregular
homeomorphisms, and
Reshetnyak was the first to ask what
information inequality \eqref{qrdef1} harbours {\it per se}.
In a series of papers in 1966--69, Reshetnyak laid the analytic
foundations for the theory of quasiregular mappings. The
single deepest fact he discovered
 was that quasiregular mappings are
branched coverings (as defined above).
  It is instructive to outline the main steps in the
proof for this remarkable assertion, which akin to Riemann's result
exerts   significant topological information from purely
analytic data. For the details,
see, {\it e.g.}, [28], [29], [18].

To wit, let $f:\Omega\to\rn$ be $K-$quasiregular.
 Fix $y\in \R^n$ and consider
the preimage $Z=f^{-1}(y)$. One
first shows that the function
$u(x)=\log|f(x)-y|$
solves a quasilinear elliptic partial differential equation
\begin{equation}\label{funny}
-{\text{div}}\mathcal{A}(x,\nabla u(x))=0, \qquad \mathcal{A}(x,\xi)\cdot \xi
\simeq |\xi|^n,
\end{equation}
in the open set $\Omega\setminus Z$ in the weak (distributional)  sense.
 In general, $\mathcal{A}$ in \eqref{funny}
depends on $f$, but
its ellipticity only on $K$ and $n$.
For holomorphic functions, {\it i.e.}, for $n=2$ and $K=1$,
equation \eqref{funny} reduces to the Laplace equation $-$div$\nabla u=0$.

Now
$u(x)$ tends to $-\infty$ continuously as $x$ tends to $Z$.
Reshetnyak develops sufficient {\textit{nonlinear potential theory}}
to conclude that such {\textit{polar sets}}, associated with  equation
\eqref{funny}, have Hausdorff dimension zero. It follows that $Z$
is totally disconnected, {\it i.e.}, the mapping $f$ is
{\textit{light}}. This is the purely analytic part of the proof. The next
step is to show that nonconstant quasiregular mappings are
{\textit{sense-preserving}}.
 This part of the proof mixes analysis and topology. What
remains is a purely topological fact that sense-preserving and light
mappings between connected oriented manifolds
are branched coverings.

Initially, Reshetnyak's theorem served as the  basis
for a higher dimensional function theory. In the 1980's, it was discovered
by researchers in nonlinear elasticity.
In the following, we shall
discuss more recent, different types of applications.

\section{The branch set}

\vskip -5mm \hspace{5mm}

Branched coverings between surfaces behave locally like analytic functions
according to a classical theorem of Sto\"\i low.
By a theorem of Chernavski\u\i, for every $n\ge2$,
the branch set of a discrete
and open mapping between $n$-manifolds has topological dimension at most
$n-2$. For branched coverings between 3-manifolds,  the branch
set is either empty or has topological dimension 1 [24], but
  in dimensions $n\ge5$ there are branched coverings between
$n$-manifolds with branch set of dimension $n-4$,
cf. Section 7.\footnote{See [23] for a recent survey on dimension theory
and branched coverings.}

The branch set of a quasiregular mapping is a somewhat enigmatic
object in dimensions $n\ge3$.
It can be very complicated, containing for example
 many wild Cantor sets
of classical geometric topology [14], [15].
There is currently no  theory available
that would  explain or describe the geometry
of allowable branch sets, cf. Problems \textbf{2} and \textbf{4} in Section 7.

In the next three sections, we shall discuss the problem of finding
bi-Lipschitz parametrizations for metric spaces. It will become
clear only later how this problem is related to
the branch set.

\section{Bi-Lipschitz parametrization of spaces}

\vskip -5mm \hspace{5mm}

A homeomorphism $f:X\to Y$ between metric spaces
is {\textit{bi-Lipschitz}} if there exists a constant $L\ge1$ such that
$$
L^{-1}d_X(a,b)\le d_Y(f(a),f(b))\le L d_X(a,b)
$$
for each pair of points $a,b\in X$.
It appears to be a difficult
problem to decide when a given a metric space $X$ can be covered
by open sets each of which is  bi-Lipschitz
homeomorphic to an open set in $\rn$,  $n\ge2$.
If this is the case, let us say,
for brevity and with a slight abuse of
language, that
$X$ is \textit{locally bi-Lipschitz equivalent to} $\rn$.

Now a separable metrizable space
 is a Lipschitz manifold (in the sense of charts)
if and only if
it admits a metric, compatible with
the given topology, that
makes the space  locally bi-Lipschitz equivalent to $\rn$ [22].
The problem here is different from characterizing Lipschitz
manifolds among topological spaces, for the metric is
given first, cf. [8], [39], [40], [41].

 To get a grasp of the difficulty of the problem,
consider the following example: {\it There exist
 finite $5$-dimensional polyhedra  that are
homeomorphic to the standard $5$-sphere ${\mathbf S^5}$, but not
locally bi-Lipschitz equivalent to $\mathbf R^5$.} This observation of
Siebenmann and Sullivan [38] is based on a deep result of Edwards
 [9], which asserts that the double suspension
$\Sigma^2H^3$ of a $3$-dimensional homology sphere
$H^3$, with nontrivial fundamental group,
is homeomorphic to the standard sphere $\mathbf S^{5}$. (See also [6].)
 One can think of
$X=\Sigma^2H^3$ as a
 join $X=\mathbf S^1*H^3$, and it is easy to check that the complement of the
{\textit{suspension circle}} $\mathbf S^1$ in $X$ is not simply connected.
Consequently, every homeomorphism
$f:X\to \mathbf S^{5}$ must transfer $\mathbf S^1$ to a closed
curve $\Gamma=f(\mathbf S^1)$ whose complement in $\mathbf S^{5}$ is
not simply connected. A
 general position
argument and  Fubini's theorem  imply that, in this case, the Hausdorff
dimension of
$\Gamma$ must be at least 3. Hence $f$ cannot be Lipschitz. In fact, $f$
cannot be H\"older continuous with any exponent greater than $1/3$.
 It is not known what other obstructions there are for a
homeomorphism $X\to \mathbf S^5$,
cf. [16, Questions 12--14].

See [33] and [37] for surveys on
parametrization and related topics.

\section{Necessary conditions}\label{problemkaksi}

\vskip -5mm \hspace{5mm}

What  are the obvious necessary conditions that a given metric
space $X$ must satisfy, if it were to be locally bi-Lipschitz
equivalent to $\R^n$, $n\ge2$? Clearly, $X$ must be an
$n$-manifold. Next, bi-Lipschitz mappings preserve Hausdorff
measure in a quantitative manner, so in particular $X$ must be
$n$-{\textit {rectifiable}} in the sense of geometric measure
theory; moreover,
 locally the Hausdorff $n$-measure
should assign to each ball of radius $r>0$ in $X$ a mass comparable
to $r^n$. Let us say that $X$ is {\emph{metrically
$n$-dimensional}} if it satisfies these geometric measure theoretic
requirements.

It is not difficult to find examples of metrically $n$-dimensional
manifolds
that are not locally  bi-Lipschitz equivalent to $\R^n$. The
measure theory allows for cusps and folds that are not tolerated
by bi-Lipschitz parametrizations. Further geometric constraints are
necessary; but, unlike in the case of
the measure theoretic conditions, it is
not
obvious what these constraints should be. A convenient
choice is that of {\textit{local linear contractibility}}:
locally each metric ball in
$X$ can be contracted to a point inside a ball with the same
center but radius multiplied by  a fixed factor.
\footnote{See [36] for  analytic implications of this condition.}

Still, a metrically
$n$-dimensional and locally linearly contractible metric $n$-manifold
need not be
locally bi-Lipschitz equivalent to  $\R^n$. The double suspension
of a  homology 3-sphere with nontrivial fundamental group
 as described in the previous
section serves as a counterexample. In 1996, Semmes
[34], [35]
exhibited examples to the same effect in all dimensions $n\ge3$, and
recently Laakso [21] crushed the last hope that the above conditions
might characterize at least
$2$-dimensional metric manifolds that are locally bi-Lipschitz
equivalent to $\R^2$.
However, unlike the examples of Edwards and Semmes,
Laakso's metric space cannot be embedded bi-Lipschitzly in any finite
dimensional Euclidean space. Thus the
following problem remains open:

\begin{problem} {\rm Let $X$ be a topological surface inside some
$\R^N$ with the inherited metric. Assume that $X$ is metrically
2-dimensional and locally linearly
 contractible. Is $X$ then locally by-Lipschitz equivalent to
 $\R^2$?}
\end{problem}

In conclusion, perhaps excepting the dimension $n=2$,
more necessary conditions are needed in order to
characterize the spaces that are
locally bi-Lipschitz equivalent
to $\rn$.\footnote{There are interesting and nontrivial
sufficient conditions known, but
these are far from being necessary [42], [43], [2], [3], [5].}
The
idea to use Reshetnyak's theorem in this connection
 originates in two papers by Sullivan [40],
[41], and is later developed in [17].
 Recall  that in this theorem topological
conclusions are drawn from purely analytic data. Now imagine that such
data would make sense in a space that is not
 {\it a priori}
Euclidean. Then, if one could obtain a branched
covering mapping into $\R^n$, {\textit{manifold points would appear, at
least outside the branch set}}.
We discuss the possibility
to develop this idea in the next section.

\section{Cartan-Whitney presentations}\label{branched}\setzero

\vskip -5mm \hspace{5mm}

Let $X$ be a metrically $n$-dimensional, linearly locally
contractible $n$-manifold that is also a metric subspace of some
$\R^N$. Suppose that there exists a bi-Lipschitz homeomorphism
$f:X\to f(X)\subset\R^n$. Then $f$  pulls back to $X$ the standard
coframe
 of $\rn$, providing almost
everywhere defined (essentially) bounded  differential 1-forms
$
\rho_i=f^*dx_i$, $i=1,\dots,n$.
To be more precise here,
by Kirzsbraun's theorem,  $f$ can be
extended to a Lipschitz mapping $\bar f:\R^N\to \R^n$, and the
1-forms
\begin{equation}\label{pull}
\rho_i=\bar f^*dx_i=d\bar f_i, \ \ \ \ \ i=1,\dots,n,
\end{equation}
are  well defined in $\R^N$ as flat 1-forms of Whitney. {\it Flat forms}
are forms with
 $L^{\infty}$-coefficients  such that the distributional exterior
differential of the form also has  $L^{\infty}$-coefficients.
The forms in \eqref{pull} are  closed, because the fundamental relation
$
d\bar f^*=\bar f^*d
$
holds true for Lipschitz maps.

 According
to a theorem of Whitney [45, Chapter IX],  flat
 forms $(\rho_i)$ have a well defined trace on $X$, and on the measurable
tangent bundle of $X$, essentially because of the
rectifiability.\footnote{There is a technical point
about orientation which we ignore here [17, 3.26].}
Because  $f=\bar f|X$ has a Lipschitz inverse,
there exists a constant $c>0$ such that
\begin{equation}\label{ala}
*(\rho_1\wedge\dots\wedge\rho_n)\ge c>0
\end{equation}
almost everywhere on $X$, where the Hodge
star operator $*$ is determined
by the chosen orientation on $X$.

Condition \eqref{ala} was turned into a definition in [17]. We
 say that $X$ admits {\textit{local
Cartan-Whitney presentations}} if for each point  $p\in X$ one can
find an
$n$-tuple of  flat
1-forms
$\rho=(\rho_1,\dots,\rho_n)$
defined in an $\R^N$-neighborhood of  $p$  such that
condition \eqref{ala} is satisfied on $X$ near the point $p$.

\begin{theorem}\label{hs}\textnormal{[17]}
Let $X\subset \R^N$ be a metrically $n$-dimensional, linearly
locally contractible $n$-manifold admitting local Cartan-Whitney
presentations.  Then
$X$ is locally bi-Lipschitz equivalent to $\R^n$ outside a
closed set of measure zero and of topological dimension at most $n-2$.
\end{theorem}

To prove Theorem \ref{hs}, fix a point $p\in X$, and let
$\rho=(\rho_1,\dots,\rho_n)$ be
a
Cartan-Whitney presentation near $p$. The
requirement that  $\rho$ be flat
together with inequality
\eqref{ala} can be seen as a quasiregularity
condition for forms.\footnote{In fact,
\eqref{ala} resembles a stronger, Lipschitz version
of \eqref{qrdef1} studied in [26], [40], [15].}
We define a
 mapping
\begin{equation}\label{reshmap}
f(x)=\int_{[p,x]}\rho
\end{equation}
for $x$ sufficiently near $p$, where $[p,x]$ is the line segment in $\R^N$
from $p$ to $x$, and claim that Reshetnyak's program can be run under the
stipulated conditions on $X$.
In particular, we show that for a sufficiently
small neighborhood $U$ of $p$ in $X$, the map
$f:U\to\R^n$ given in \eqref{reshmap} is a branched covering which
is locally bi-Lipschitz outside its branch set $B_{f}$, which
furthermore is   of measure zero and of
 topological dimension
 at most $n-2$. It is important to note that $\rho$ is not assumed
to be closed, so that $df\ne\rho$ in general.

In executing Reshetnyak's proof, we use recent advances
of differential analysis
on nonsmooth spaces [13], [20], [36], as well as the theory
developed simultaneously in [15]. Incidentally,
we avoid the use of the
 Harnack
inequality for
solutions, and therefore a deeper
use of equation \eqref{funny}; this small
 improvement to Reshetnyak's argument was
found  earlier in a different context in [12].

Theorem \ref{hs} provides bi-Lipschitz coordinates for $X$ only on
a dense open set. In general, one cannot have more than that.
The double
suspension of a homology 3-sphere, as discussed
in Section 4, can be mapped to the standard 5-sphere by a
finite-to-one, piecewise linear sense-preserving map.
By pulling back the
standard coframe by such map, we obtain a global Cartan-Whitney
presentation on
a space that is not locally bi-Lipschitz equivalent to
$\mathbf R^5$.
Similar examples in dimension $n=3$
were constructed in [14], [15], by using
Semmes's spaces [34], [35].
On the other hand, we have the following result:

\begin{theorem}\label{unpub}
Let $X\subset \R^N$ be a metrically $2$-dimensional, linearly
locally contractible $2$-manifold admitting local Cartan-Whitney
presentations.  Then
$X$ is locally bi-Lipschitz equivalent to $\R^2$.
\end{theorem}

Theorem \ref{unpub} is an observation
of M. Bonk and myself.
We use
Theorem
 \ref{hs} together with the observation that, in dimension $n=2$,
the branch set consists of isolated points, which can be resolved.
The resolution follows from the
measurable Riemann mapping theorem together with the
recent work by Bonk and Kleiner [4].
 While Theorem
\ref{unpub} presents a characterization of surfaces in Euclidean space
that admit local bi-Lipschitz coordinates, we do not know whether the
stipulation about the existence of  local Cartan-Whitney
presentations is really necessary (compare Problem 1 and the
discussion preceding it).

For dimensions $n\ge3$, it would be interesting
to know when there is
 no branching in the map \eqref{reshmap}.
In [17], we ask if this be the
case when the  flat forms $(\rho_i)$ of the Cartan-Whitney
presentation
 belong to a Sobolev space $H^{1,2}_{loc}$
on $X$.  The relevant example here is the map $(r,\theta,z)\mapsto
(r,2\theta,z)$, in the cylindrical coordinates of $\rn$,
 which pulls back the standard coframe to a frame that lies
in the Sobolev space $H^{1,2-\epsilon}_{loc}$ for each $\epsilon>0$.
Indeed, it was shown in [11] that in $\rn$
every (Cartan-Whitney)
 pullback frame in $H^{1,2}_{loc}$ must come from a locally injective
mapping.


\section{Other recent progress and open problems}\label{problems}

\vskip -5mm \hspace{5mm}

In his 1978 ICM address, V\"ais\"al\"a [44] asked
 whether the branch set of a $C^1$-smooth
quasiregular mapping is empty if $n\ge3$.
It was known that $C^{n/(n-2)}$-smooth
quasiregular mappings have no branching when $n\ge3$. The proof
in [Ri, p. 12] of this fact
uses quasiregularity in a rather minimal way. In this light,
the following recent result may appear
surprising:

\begin{theorem}\label{smooth}\textnormal{[1]}
For every $\epsilon>0$ there exists a degree two
$C^{3-\epsilon}$-smooth
quasiregular mapping
$f:\cS^3\to\cS^3$ with branch set homeomorphic to $\cS^1$.
\end{theorem}

We are also able to improve the previous results as follows:

\begin{theorem}\label{tarkka}\textnormal{[1]}
Given $n\ge3$ and $K\ge1$, there exist
$\epsilon=\epsilon(n,K)>0$ and $\epsilon'=\epsilon'(n,K)>0$ such that
the branch set of every
$K$-quasiregular mapping in a domain in $\rn$ has Hausdorff dimension at
most $n-\epsilon$, and that every
$C^{n/(n-2)-\epsilon'}$-smooth
$K$-quasiregular mapping in a domain in $\rn$ is a local homeomorphism.
\end{theorem}

The second assertion in Theorem \ref{tarkka}
follows from the first, by way of
 Sard-type techniques.
The first assertion  was known earlier in a local form where $\epsilon>0$
was dependent  on the local degree [31]. Our
improvement uses [31] together with the work [30] by Rickman and Srebro.

The methods in [1] fall short in showing the sharpness
 of Theorem \ref{tarkka} in dimensions $n\ge4$
in two technical aspects.
First, we would
need to construct a quasiconformal homeomorphism of $\rn$ to itself
that is uniformly expanding  on a codimension two
affine subspace; moreover, such a map needs to be smooth outside this
subspace.
In $\R^3$, it is easier to construct a mapping with
 expanding behavior on a line; moreover, every quasiconformal
homeomorphism in dimension three can be
smoothened (with bounds) outside a given closed set [19].

We finish with some open  problems related to branching
and quasiregular mappings. The problems are
 neither new nor due to the author.

\begin{problem} {\rm What are  the
possible values for the topological dimension of the branch set
of a quasiregular mapping?}
\end{problem}

By  suspending a covering map $H^3\to \cS^3$, where $H^3$ is as in Section 4,
and using Edwards's theorem,
one finds that
there exists a branched covering
$\cS^5\to\cS^5$ that branches exactly on $\cS^1\subset \cS^5$.
It is not known whether there exists a quasiregular mapping $\cS^5\to\cS^5$
with similar branch set.
If no such map existed,  we would have an interesting implication
to a seemingly unrelated parametrization problem;  it would
follow that no double suspension of a  homology
3-sphere with nontrivial fundamental group
admits a {\it{quasisymmetric}} homeomorphism onto the standard 5-sphere,
cf. [38], [16, Question 12].

By work of Bonk and Kleiner [4],  the bi-Lipschitz
parametrization problem in dimension $n=2$
 is equivalent to an analytic
problem of characterizing, up to a bounded factor,
 the Jacobian determinants of quasiconformal
mappings in $\mathbf R^2$.
An affirmative answer to Problem
1 in Section 5 would give an affirmative answer to the following problem.

\begin{problem} \textnormal{(Compare [16, Question 2])
Is every $A_1$-weight in $\mathbf R^2$ locally comparable to the
Jacobian determinant of a quasiconformal mapping?}
\end{problem}

An $A_1$-{\it weight} is a nonnegative locally integrable function
whose mean-value over each ball is comparable to its essential
infimum over the ball.  See [7], [32],
[3], [16]  for further discussion of this and related problems.

\begin{problem}\textnormal{[16, Question 28] Is there   a
  branched
covering $f:\cS^n\to\cS^n$, for some $n\ge3$, such that
 for
every pair of homeomorphisms $\phi,\psi:\cS^n\to\cS^n$, the
mapping $\phi\circ f\circ \psi$ fails to be quasiregular?}
\end{problem}

Branched coverings constructed  by using the double suspension are
obvious candidates for such mappings.
 In [15, 9.1], we give an example
of  a branched
covering $f:\cS^3\to\cS^3$ such that for every homeomorphism
$\psi:\cS^3\to\cS^3$,
$ f\circ
\psi$  fails to be quasiregular. The example is based on a geometric
decomposition
space arising from  {\it Bing's double} [34].

We close this article by commenting on the lack of direct proofs
for some fundamental properties of quasiregular mappings related
to branching.
For example, it
is known that for each $n\ge3$ there exists $K(n)>1$ such that
every $K(n)$-quasiregular mapping is a local homeomorphism [25],
[28, p. 232].
All known proofs for this fact are indirect, exploiting the
Liouville theorem,
 and in
particular there is no numerical estimate for $K(n)$. It has been
conjectured that the winding mapping $(r,\theta, z)\mapsto
(r,2\theta, z)$  is the extremal here
(cf.  Section 6). Thus, if one
uses the {\it inner dilatation} $K_I(f)$ of a quasiregular mapping,
then conjecturally $K_I(f)<2$ implies that $B_f=\emptyset$ for a
quasiregular
 mapping $f$ in $\rn$ for $n\ge3$ [29, p. 76].

Ostensibly different, but obviously a related issue, arises in search
of {\it Bloch's constant} for quasiregular mappings. Namely,
 by exploiting normal families, Eremenko
[10]
recently proved  that for given
$n\ge3$ and $K\ge1$, there exists
$b_0=b_0(n,K)>0$ such that every $K$-quasiregular mapping $f:\rn\to\cS^n$
has an inverse branch in some ball in $\cS^n$ of radius $b_0$. No numerical
estimate for $b_0$ is known.  More generally, despite the deep results
on value
distribution of quasiregular mappings,
 uncovered by Rickman over the past quarter century,
the affect of branching on value distribution is unknown, cf. [29, p. 96].

\label{lastpage}

\end{document}